\documentclass{amsart}

\usepackage{amssymb}
\usepackage{graphicx}
\usepackage{amsmath}
\usepackage{amsthm}
\usepackage{xcolor}
\usepackage{tikz}
\usetikzlibrary{arrows}

\newcommand{\bl}[1]{\textcolor{blue}{#1}}

\theoremstyle{plain}
\newtheorem{theorem}{Theorem}[section]
\newtheorem{corollary}[theorem]{Corollary}
\newtheorem{lemma}[theorem]{Lemma}
\newtheorem{proposition}[theorem]{Proposition}

\theoremstyle{definition}
\newtheorem{definition}[theorem]{Definition}
\newtheorem{example}[theorem]{Example}
\newtheorem{remark}[theorem]{Remark}

\usepackage[hidelinks]{hyperref}
\hypersetup{colorlinks=true,citecolor=.,linkcolor=.,urlcolor=.}

\makeatletter
\@namedef{subjclassname@2020}{\textup{2020} Mathematics Subject Classification}
\makeatother

\begin{document}

\title[A note on commutative gi-lattices]
{A note on distributive lattices with a commutative generalized implication}

\author{Ismael Calomino}
\address{Facultad de Ciencias y Artes, Universidad Cat\'olica de \'Avila, C/Canteros s/n, 05005 \'Avila, Espa\~{n}a.}
\email{ismaelm.calomino@ucavila.es}

\author{Sergio Celani}
\address{CONICET and Facultad de Ciencias Exactas, Universidad Nacional del Centro, Pinto 399 Tandil, Argentina.}
\email{scelani@exa.unicen.edu.ar}


\begin{abstract}
We study distributive lattices with an ideal-valued generalized implication satisfying a commutativity identity. This class extends commutative weak Heyting algebras by allowing implication values to be non-principal ideals. We characterize its intersection with several known classes of gi-lattices, study gi-deductive systems and their relationship with gi-filters and gi-congruences, and extend the relational characterization of commutativity to this setting. The resulting correspondence yields a Priestley-style duality by restricting the established duality for gi-lattices.
\end{abstract}

\maketitle

\section{Introduction and preliminaries}

Distributive lattices with a generalized implication, or gi-lattices, were introduced in \cite{Castro2} (see also \cite{Celani-Hernan2}) and are closely related to relative annihilators \cite{Mandelker}, weak Heyting algebras \cite{Celani-Jansana}, quasi-modal algebras \cite{Castro1,Celani-QM}, as well as to more recent developments on subordination relations on bounded distributive lattices \cite{Guram,Calomino,Celani-Sub}. A gi-lattice is a structure $\langle A, \Rightarrow \rangle$, where $A$ is a bounded distributive lattice and $\Rightarrow \colon A \times A \to {\rm{Id}}(A)$ assigns an ideal $a \Rightarrow b$ to each pair $(a,b)$. Since a generalized implication takes values in $ {\rm{Id}}(A)$ rather than in $A$, it is not a binary operation on $A$. Thus gi-lattices are not algebras of the usual single-sorted type, although they share many structural properties with weak Heyting algebras and quasi-modal lattices.

A weak Heyting algebra, or WH-algebra, is defined as a pair $\langle A, \to \rangle$, where $A$ is a bounded distributive lattice and $\to \colon A \times A \to A$ is a function such that the following conditions are satisfied for all $a,b,c \in A$:
\begin{enumerate}
\item $ a \to (b \wedge c) = (a \to b) \wedge (a \to c)$,
\item $(a \vee b) \to c = (a \to c) \wedge (b \to c)$,
\item $(a \to b) \wedge (b \to c) \leq a \to c$,
\item $a \to a = 1$.
\end{enumerate}
The variety of WH-algebras is the algebraic counterpart of the least subintuitionistic logic $wK$ \cite{Celani-Jansana}. Moreover, if a WH-algebra $\langle A, \to \rangle$ additionally satisfies the inequalities $x \leq 1 \to x$, $x \to y \leq z \to (x \to y)$ and $x \wedge (x \to y) \leq y$, then it is a Heyting algebra. Examples of WH-algebras include the basic algebras introduced in \cite{Ardeshir-Ruitenburg} (see also \cite{Alizadeh-Ardeshir}) and the subresiduated lattices introduced in \cite{Epstein-Horn}. Every weak Heyting algebra induces a gi-lattice structure in the following way: the binary map from $A \times A$ to the set of all ideals of $A$ defined by $a \Rightarrow b = (a \rightarrow b]$ is a generalized implication on $A$, making the structure $\langle A, \Rightarrow \rangle$ a gi-lattice. Nourany, Ghorbani and Borumand Saeid introduced self-distributive weak Heyting algebras in \cite{Ghorbani2} and studied their filters and congruences in \cite{Ghorbani1}. In these algebras the strict implication satisfies $a\to(b\to c)=(a\to b)\to(a\to c)$. This identity is equivalent to $a\to(b\to c)=b\to(a\to c)$ \cite[Proposition 3.8]{Ghorbani2}, which explains the alternative name commutative weak Heyting algebras. Celani and San Mart\'in established a Priestley-style duality for commutative weak Heyting algebras in \cite{Celani-Hernan}, building on the duality for weak Heyting algebras in \cite{Celani-Jansana}. Our aim is to extend this study to ideal-valued generalized implications, that is, to gi-lattices that satisfy the equation $a \Rightarrow (b \Rightarrow c) = b \Rightarrow (a \Rightarrow c)$. The distinction matters for infinite lattices: every ideal of a finite lattice is principal, whereas an implication value in an infinite gi-lattice need not have a greatest element. The relational conditions studied below extend those in \cite{Celani-Hernan}; the filter-theoretic results extend, and should be compared with, the corresponding results of \cite{Ghorbani1}. The underlying notions of gi-lattice, gi-filter and gi-congruence, and the general duality, are those of \cite{Castro2}.

On the other hand, generalized implication provides a direct connection to non-classical logic. Makinson and van der Torre introduced input/output logic in \cite{Makinson-Torre} as a formal framework for modeling the interaction between logical inference and other agent-related notions, with its primary application being in the context of formalizing normative systems. In input/output logic, a system consists of relational pairs where an input condition $a$ generates an output obligation $b$. The study of input/output logic from an algebraic perspective is grounded in subordination algebras in \cite{Domenico1,Domenico2}, whose normative properties correspond to generalized implications and are closely related to slanted algebras \cite[Example 3.3]{DeRudder-Palmigiano}. Generalized implication takes the input and, instead of returning a single value, generates a set, meaning that the valid output is no longer an isolated element, but rather all possible logical consequences grouped into an ideal. Furthermore, adding the commutativity condition to generalized implication means that the order in which contexts are processed does not alter the final obligations. This algebraic property guarantees that if a norm depends on two independent conditions, the final normative outcome will be the same regardless of which one is evaluated first, ensuring that different deductive paths converge to the same ideal of outputs.

The paper is organized as follows. In the rest of this section, we recall the basic notions and results that will be used throughout the paper. In Section \ref{sec2}, we recall the definition of gi-lattice and establish some results needed in the following sections. In Section \ref{sec3}, we introduce the class of commutative gi-lattices, study their relationship with other subclasses of gi-lattices, and provide some characterizations. In Section \ref{sec4}, we present the notion of gi-deductive system and its connection with gi-congruences. Finally, Section \ref{sec5} obtains a Priestley-style duality by restricting the duality for gi-lattices to the commutative objects.

Let $\langle X, \leq \rangle$ be a poset. We say that $U \subseteq X$ is an upset if for every $x,y \in X$, if $x \leq y$ and $x \in U$, then $y \in U$. Analogously, we say that $U \subseteq X$ is a downset if for every $x,y \in X$, if $x \leq y$ and $y \in U$, then $x \in U$. We write {\rm{Up}}$(X)$ for the poset of upsets of $X$. If $U \subseteq X$, then we write $[U)$ for the upset generated by $U$ and $(U]$ for the downset generated by $U$. Thus, $[U) = \{ x \in X \colon \exists y \in U (y \leq x) \}$ and $(U] = \{ x \in X \colon \exists y \in U (x \leq y) \}$. If $U = \{u\}$ we write $[u)$ in place of $[\{u\})$ and $(u]$ in place of $(\{u\}]$. It follows that $[u) = \{ y \in X \colon u \leq y \}$ and $(u] = \{ y \in X \colon y \leq u \}$.

Let $A$ be a bounded lattice. We say that $F \subseteq A$ is a filter if $1 \in F$, $F$ is an upset, and $a \wedge b \in F$ whenever $a,b \in F$. Dually, we say that $I \subseteq A$ is an ideal if $0 \in I$, $I$ is a downset, and $a \vee b \in I$ whenever $a,b \in I$. We denote by ${\rm{Fi}}(A)$ and ${\rm{Id}}(A)$ the sets of all filters and ideals of $A$, respectively. For every $X \subseteq A$ we write ${\rm{Fig}}(X)$ to indicate the filter generated by $X$. It is also known that the filter generated by a non-empty subset $X$ of $A$ can be characterized by
\[
{\rm{Fig}}(X) = \{ a \in A \colon \exists x_1, \dots, x_n \in X (x_1 \wedge \dots \wedge x_n \leq a)\}.
\]
Analogously, the ideal ${\rm{Idg}}(X)$ generated by a non-empty subset $X$ of $A$ can be characterized by
\[
{\rm{Idg}}(X) = \{ a \in A \colon \exists x_1, \dots, x_n \in X (a \leq x_1 \vee \dots \vee x_n)\}.
\]
Note that, for every $a \in A$, ${\rm{Fig}}(\{a\}) = [a)$ and ${\rm{Idg}}(\{a\}) = (a]$. A filter $P$ of $A$ is said to be prime if it is proper and for every $a,b \in A$ it holds that if $a \vee b \in P$, then $a \in P$ or $b \in P$. We denote by $X(A)$ the set of all prime filters of $A$.

We recall Priestley duality \cite{Priestley,Balbes}. A Priestley space is a triple $\langle X, \leq, \tau \rangle$ such that $\langle X, \leq \rangle$ is a poset, $\langle X, \tau \rangle$ is a compact topological space and the Priestley separation axiom is satisfied, which means that for every $x,y \in X$, if $x \nleq y$, then there exists a clopen upset $U$ such that $x \in U$ and $y \notin U$. The family of clopen upsets of $X$ is denoted by $CU(X)$, and it is well known that $CU(X)$ is a bounded distributive lattice. Consider the map $\varphi \colon A \to {\rm{Up}}(X(A))$ given by $\varphi(a) = \{ P \in X(A) \colon a \in P \}$, for each $a \in A$. Then the triple $\langle X(A), \subseteq, \tau_{A} \rangle$ is a Priestley space, where $\tau_{A}$ is the topology generated by the subbase $\{ \varphi(a) \colon a \in A \} \cup \{ \varphi(a)^{c} \colon a \in A \}$. Moreover, since $CU(X(A)) = \{ \varphi(a) \colon a \in A \}$, the lattices $A$ and $CU(X(A))$ are isomorphic. If $X$ is a Priestley space, then the map $\epsilon \colon X \to X(CU(X))$ given by $\epsilon(x) = \{ U \in CU(X) \colon x \in U \}$ is a homeomorphism and an order-isomorphism. If $A$ and $B$ are bounded distributive lattices and $h \colon A \to B$ is a homomorphism, then the map $h^{\ast} \colon X(B) \to X(A)$ defined by $h^{\ast}(P) = h^{-1}(P)$ is a continuous and monotone map. Conversely, if $X$ and $Y$ are Priestley spaces and $f \colon X \to Y$ is a continuous and monotone map, then the map $f^{\ast} \colon CU(Y) \to CU(X)$ defined by $f^{\ast}(U) = f^{-1}(U)$ is a homomorphism of bounded distributive lattices. Hence, there is a dual equivalence between the category of bounded distributive lattices with homomorphisms and the category of Priestley spaces with continuous and monotone maps.

\section{Generalized implication} \label{sec2}

In this section, we recall generalized implications on bounded distributive lattices and establish some auxiliary results.

\begin{definition} \label{def: generalized implication}
Let $A$ be a bounded distributive lattice. A generalized implication on $A$ is a function $\Rightarrow \colon A \times A \to {\rm{Id}}(A)$ that satisfies the following conditions for all $a,b,c \in A$:
\begin{enumerate}
\item[(GI1)] $(a \Rightarrow b) \cap (a \Rightarrow c) = a \Rightarrow (b \wedge c)$,
\item[(GI2)] $(a \Rightarrow c) \cap (b \Rightarrow c) = (a \vee b) \Rightarrow c$,
\item[(GI3)] $(a \Rightarrow b) \cap (b \Rightarrow c) \subseteq a \Rightarrow c$,
\item[(GI4)] $a \Rightarrow a = A$.
\end{enumerate}
A distributive lattice with a generalized implication, or gi-lattice for short, is a pair $\langle A, \Rightarrow \rangle$, where $A$ is a bounded distributive lattice and $\Rightarrow$ is a generalized implication defined on $A$.
\end{definition}

We denote by $\bf{GI}$ the class of gi-lattices. We recall several examples and constructions.

\begin{example} \label{Example_1}
Let $A$ be a bounded distributive lattice. A quasi-modal operator on $A$ is a function $\Delta \colon A \to {\rm{Id}}(A)$ such that $\Delta 1 = A$ and $\Delta (a \wedge b) = \Delta a \cap \Delta b$, for every $a,b \in A$ \cite{Castro1,CelaniMonteiro,Celani-QM}. A quasi-modal lattice is a structure $\langle A, \Delta \rangle$, where $A$ is a bounded distributive lattice and $\Delta$ is a quasi-modal operator on $A$. If $\langle A, \Rightarrow \rangle$ is a gi-lattice, then the function $\Delta \colon A \to {\rm{Id}}(A)$ given by $\Delta a = 1 \Rightarrow a$ is a quasi-modal operator on $A$. Thus, the pair $\langle A, \Delta \rangle$ is a quasi-modal lattice.
\end{example}

\begin{example} \label{Example_3}
Let $A$ be a lattice and $a,b \in A$. The annihilator of $a$ relative to $b$ is the set $\langle a, b \rangle = \{ x \in A \colon x \wedge a \leq b \}$. A well-known result due to Mandelker \cite{Mandelker} asserts that $A$ is distributive if and only if $\langle a, b \rangle$ is an ideal of $A$, for every $a,b \in A$. Hence, if $A$ is a bounded distributive lattice, then the map $\Rightarrow \colon A \times A \to {\rm{Id}}(A)$ given by $a \Rightarrow b = \langle a, b \rangle$ is a generalized implication on $A$.
\end{example}

\begin{example} \label{Example_2}
If $\langle A, \to \rangle$ is a WH-algebra, then the function $\Rightarrow \colon A \times A \to {\rm{Id}}(A)$ defined by $a \Rightarrow b = (a \to b]$ is a generalized implication on $A$.
\end{example}

\begin{example} \label{Example_conj}
Let $\langle X, \leq, R \rangle$ be a structure such that $\langle X, \leq \rangle$ is a poset and $R \subseteq X \times X$ is a binary relation satisfying $\leq \circ R \subseteq R$. If we define $\to_{R} \colon {\rm{Up}}(X) \times {\rm{Up}}(X) \to {\rm{Up}}(X)$ by
\begin{equation} \tag{$\bullet$}
U \to_{R} V = \{ x \in X \colon R(x) \cap U \subseteq V \}, \label{impli conj}
\end{equation}
for every $U,V \in {\rm{Up}}(X)$, then $\langle {\rm{Up}}(X), \cap, \cup, \to_{R}, \emptyset, X \rangle$ is a WH-algebra. So, by Example \ref{Example_2}, the pair $\langle {\rm{Up}}(X), \Rightarrow_{R} \rangle$ is a gi-lattice where $U \Rightarrow_{R} V = ( U \to_{R} V]$, for every $U,V \in {\rm{Up}}(X)$.
\end{example}

\begin{example} \label{Example_constant}
Let $A$ be a bounded distributive lattice and let $I \in {\rm{Id}}(A)$. Define $\Rightarrow_{I} \colon A \times A \to {\rm{Id}}(A)$ by
\[
a \Rightarrow_{I} b =
\begin{cases}
A & \text{if } a \leq b,\\
I & \text{if } a \nleq b.
\end{cases}
\]
Then $\Rightarrow_{I}$ is a generalized implication on $A$. Indeed, all the values are ideals and (GI4) is immediate. For (GI1), note that $a \leq b \wedge c$ if and only if $a \leq b$ and $a \leq c$; hence $a \Rightarrow_{I} (b \wedge c) = A$ if and only if $(a \Rightarrow_{I} b) \cap (a \Rightarrow_{I} c) = A$, and both sides are equal to $I$ otherwise. Condition (GI2) follows in the same way from the fact that $a \vee b \leq c$ if and only if $a \leq c$ and $b \leq c$. For (GI3), if $a \leq c$ there is nothing to prove; if $a \nleq c$, then $a \leq b$ and $b \leq c$ cannot both hold, so at least one of $a \Rightarrow_{I} b$, $b \Rightarrow_{I} c$ is equal to $I$, and therefore $(a \Rightarrow_{I} b) \cap (b \Rightarrow_{I} c) \subseteq I = a \Rightarrow_{I} c$. When $I = \{0\}$ all the values are principal, and $\Rightarrow_{\{0\}}$ is the implication associated with the WH-algebra $\langle A, \to \rangle$ given by $a \to b = 1$ if $a \leq b$ and $a \to b = 0$ otherwise; the gi-lattice of Remark \ref{remark_CGI} below is of this form. If $I$ is not principal, then $a \Rightarrow_{I} b$ is not principal whenever $a \nleq b$, so $\langle A, \Rightarrow_{I} \rangle$ does not arise from any WH-algebra as in Example \ref{Example_2}. For instance, let $A = {\rm{FC}}(\mathbb{N})$ be the Boolean algebra of all finite and cofinite subsets of $\mathbb{N}$ and let $I = {\rm{Fin}}$ be the ideal of all finite subsets of $\mathbb{N}$, which has no greatest element. Then $X \Rightarrow_{{\rm{Fin}}} Y = {\rm{Fin}}$ whenever $X \nsubseteq Y$. We will see in Remark \ref{remark_constant_not_CGI} that these gi-lattices are, in general, not commutative.
\end{example}

\begin{proposition} \cite{Castro2} \label{prop: properties}
Let $\langle A, \Rightarrow \rangle \in \bf{GI}$ and $a,b,c \in A$. If $a \leq b$, then $(c \Rightarrow a) \subseteq (c \Rightarrow b)$, $(b \Rightarrow c) \subseteq (a \Rightarrow c)$ and $a \Rightarrow b = A$. Moreover, $(a \Rightarrow b) \cap (a \Rightarrow c) \subseteq a \Rightarrow (b \vee c)$.
\end{proposition}

Let $\langle A, \Rightarrow \rangle \in \bf{GI}$. We consider the binary relation $R_{A} \subseteq {\rm{Fi}}(A) \times {\rm{Fi}}(A)$ given by
\begin{equation} \tag{$\star$}
(F,G) \in R_{A} \Longleftrightarrow \forall a,b \in A \left( (a \Rightarrow b) \cap F \neq \emptyset \text{ and } a \in G \text{ imply } b \in G \right). \label{relation R_A}
\end{equation}
Let $F \in {\rm{Fi}}(A)$. We define the operator $D_{F} \colon {\mathcal{P}}(A) \to {\mathcal{P}}(A)$ by
\[
D_{F}(X) = \left\{ b \in A \colon \exists Y \subseteq X \text{ finite such that }  \left( \bigwedge Y \Rightarrow b \right) \cap F \neq \emptyset \right\}.
\]
Here $\bigwedge\emptyset=1$. The map $D_F$ is a finitary closure operator, and $D_F(X)$ is the least filter containing $X$ such that $(F,D_F(X))\in R_A$ \cite[Proposition 11]{Castro2}. In particular, $D_{F}(X) \in {\rm{Fi}}(A)$ and $\left( F, D_{F}(X) \right) \in R_{A}$, for every  $X \subseteq A$. The following result is a generalization of the prime filter theorem for bounded distributive lattices.

\begin{theorem} \cite{Castro2} \label{theorem_FP}
Let $\langle A, \Rightarrow \rangle \in \bf{GI}$. Let $X \subseteq A$, $F \in {\rm{Fi}}(A)$ and $I \in {\rm{Id}}(A)$. If $D_{F}(X) \cap I = \emptyset$, then there exists $P \in X(A)$ such that $D_{F}(X) \subseteq P$, $(F,P) \in R_{A}$ and $P \cap I = \emptyset$.
\end{theorem}

\begin{corollary} \cite{Castro2} \label{corollary_FP}
Let $\langle A, \Rightarrow \rangle \in \bf{GI}$. Let $a,b \in A$ and $P \in X(A)$. Then $(a \Rightarrow b) \cap P = \emptyset$ if and only if there exists $Q \in X(A)$ such that $(P,Q) \in R_{A}$, $a \in Q$ and $b \notin Q$.
\end{corollary}

Let $\langle A, \Rightarrow \rangle \in {\bf{GI}}$. Let $F \in {\rm{Fi}}(A)$ and $I \in {\rm{Id}}(A)$. We define the following ideal:
\[
F \Rightarrow I = {\rm{Idg}} \left( \bigcup \{ f \Rightarrow i \colon f \in F \text{ and } i \in I \} \right).
\]

\begin{lemma} \label{lema_F to I} \cite{Calomino}
Let $\langle A, \Rightarrow \rangle \in {\bf{GI}}$. Let $F ,G \in {\rm{Fi}}(A)$, $I \in {\rm{Id}}(A)$ and $a,b \in A$. Then:
\begin{enumerate}
\item $F \Rightarrow I = \{ x \in A \colon \exists f \in F, \, \exists i \in I \, (x \in f \Rightarrow i) \}$.
\item $[a) \Rightarrow (b] = a \Rightarrow b$.
\item $D_{G}(F) \cap I = \emptyset$ if and only if $G \cap (F \Rightarrow I) = \emptyset$.
\item $G \cap (F \Rightarrow I) = \emptyset$ if and only if there exists $P \in X(A)$ such that $F \subseteq P$, $(G,P) \in R_{A}$ and $P \cap I = \emptyset$.
\end{enumerate}
\end{lemma}

In view of Lemma \ref{lema_F to I}(2), for $a \in A$ and $I \in {\rm{Id}}(A)$ we write $a \Rightarrow I$ instead of $[a) \Rightarrow I$. In particular, expressions such as $a \Rightarrow (b \Rightarrow c)$ or $z \Rightarrow (x \Rightarrow y)$ denote ideals of $A$. Moreover, if $a \leq f$, then $f \Rightarrow i \subseteq a \Rightarrow i$ by Proposition \ref{prop: properties}, so Lemma \ref{lema_F to I}(1) yields
\begin{equation} \tag{$\ast$}
a \Rightarrow I = \{ x \in A \colon \exists i \in I \, (x \in a \Rightarrow i) \} = \bigcup \{ a \Rightarrow i \colon i \in I \}. \label{a to I}
\end{equation}
Consequently, if $a \leq b$ and $I \subseteq J$, then $b \Rightarrow I \subseteq a \Rightarrow J$.

\subsection{The classes ${\bf{TGI}}$, ${\bf{RGI}}$, ${\bf{TRGI}}$ and ${\bf{BasGI}}$}

We recall some particular classes of gi-lattices that were studied in \cite{Calomino}:
\[
\begin{array}{rcccl}
{\bf{TGI}}            & = &  {\bf{GI}}  & + &  \{ x \Rightarrow y \subseteq z \Rightarrow ( x \Rightarrow y) \}, \\
{\bf{RGI}}            & = &  {\bf{GI}}  & + &  \{ (x] \cap (x \Rightarrow y) \subseteq (y] \}, \\
{\bf{TRGI}}          & = &  {\bf{GI}}  & + &  \{ x \Rightarrow y \subseteq z \Rightarrow ( x \Rightarrow y), (x] \cap (x \Rightarrow y) \subseteq (y] \},  \\
{\bf{BasGI}}         & = &  {\bf{GI}}  & + &  \{ (x] \subseteq 1 \Rightarrow x \}. \\
\end{array}
\]
The class ${\bf{BasGI}}$ is a generalization of the basic algebras introduced by Ardeshir and Ruitenburg in \cite{Ardeshir-Ruitenburg} and the class ${\bf{TRGI}}$ is a generalization of the subresiduated lattices introduced by Epstein and Horn in \cite{Epstein-Horn}. Also, the classes ${\bf{TGI}}$ and ${\bf{RGI}}$ generalize the TWH-algebras and RWH-algebras, respectively, studied in \cite{Celani-Jansana}.

\begin{theorem} \label{characterization R_A} \cite{Calomino}
Let $\langle A, \Rightarrow \rangle \in {\bf{GI}}$. Then:
\begin{enumerate}
\item $\langle A, \Rightarrow \rangle \in {\bf{TGI}}$ if and only if $R_{A}$ is transitive.
\item $\langle A, \Rightarrow \rangle \in {\bf{RGI}}$ if and only if $R_{A}$ is reflexive.
\item $\langle A, \Rightarrow \rangle \in {\bf{BasGI}}$ if and only if for all $F,G \in {\rm{Fi}}(A)$, if $(F,G) \in R_{A}$ then $F \subseteq G$.
\end{enumerate}
\end{theorem}

\begin{corollary} \label{DL} \cite{Calomino}
Let $\langle A, \Rightarrow \rangle \in {\bf{GI}}$. Then:
\begin{enumerate}
\item $\langle A, \Rightarrow \rangle \in {\bf{RGI}}$ if and only if $x \in y \Rightarrow z$ implies $x \wedge y \leq z$.
\item $\langle A, \Rightarrow \rangle \in {\bf{BasGI}}$ if and only if $x \wedge y \leq z$ implies $x \in y \Rightarrow z$.
\end{enumerate}
\end{corollary}

\begin{remark} \label{remark_DL0}
It follows that ${\bf{BasGI}} \subseteq {\bf{TGI}}$. On the other hand, if $\langle A, \Rightarrow \rangle \in {\bf{BasGI}} \cap {\bf{RGI}}$, then by Corollary \ref{DL} we have that $a \wedge b \leq c$ if and only if $a \in b \Rightarrow c$, i.e., the generalized implication $\Rightarrow$ is the relative annihilator. We have the following diagram of inclusions among these subclasses of ${\bf{GI}}$:
\vspace{0.4cm}
\begin{center}
\begin{tikzpicture}
\filldraw [thick]
(0,1) circle (2pt) node[left] {${\bf{BasGI}}$} --
(1,0) circle (2pt) node[right] {${\bf{DL}}$}     --
(1,0)--(2,1) circle (2pt) node[right] {${\bf{TRGI}}$} --
(1,2) circle (2pt) node[left] {${\bf{TGI}}$}--(0,1)
(2,1)--(3,2)circle (2pt) node[right] {${\bf{RGI}}$}--
(2,3)circle (2pt) node[left]{${\bf{GI}}$} -- (1,2);
\end{tikzpicture}
\end{center}
\vspace{0.2cm}
where ${\bf{DL}}$ denotes the class of gi-lattices whose generalized implication is the relative annihilator (which can be identified with the class of bounded distributive lattices). For more details see \cite{Calomino}.
\end{remark}

\begin{proposition} \label{propo_annihilator}
Let $\langle A, \Rightarrow \rangle \in {\bf{GI}}$. Then the generalized implication $\Rightarrow$ is the relative annihilator if and only if $1 \Rightarrow a = (a]$, for every $a \in A$.
\end{proposition}
\begin{proof}
If $a \Rightarrow b = \langle a, b \rangle$, then $1 \Rightarrow b = \langle 1, b \rangle = \{ x \in A \colon x \leq b \} = (b]$. Conversely, suppose that $1 \Rightarrow a = (a]$, for every $a \in A$. We show that $\langle A, \Rightarrow \rangle \in {\bf{BasGI}} \cap {\bf{RGI}}$. Since $1 \Rightarrow a = (a]$, we have $(a] \subseteq 1 \Rightarrow a$ and so $\langle A, \Rightarrow \rangle \in {\bf{BasGI}}$. On the other hand, by hypothesis and (GI3) of Definition \ref{def: generalized implication},
\[
(a] \cap (a \Rightarrow b) = (1 \Rightarrow a) \cap (a \Rightarrow b) \subseteq 1 \Rightarrow b = (b].
\]
Then $\langle A, \Rightarrow \rangle \in {\bf{RGI}}$ and by Remark \ref{remark_DL0}, the generalized implication $\Rightarrow$ is the relative annihilator.
\end{proof}

The following auxiliary results will be used in Section \ref{sec3}. The first extends the corresponding prime-filter statement for weak Heyting algebras in \cite{Celani-Hernan}; the second is already contained in \cite[Proposition 11(5)]{Castro2}.

\begin{lemma} \label{lema_aux1}
Let $\langle A, \Rightarrow \rangle \in {\bf{GI}}$ and $P,Q,R \in X(A)$. Then the following conditions are equivalent:
\begin{enumerate}
\item $D_{P}(Q) \subseteq R$.
\item There exists $Z \in X(A)$ such that $Q \subseteq Z \subseteq R$ and $(P,Z) \in R_{A}$.
\end{enumerate}
\end{lemma}
\begin{proof}
$(1)\Rightarrow(2)$: Since $R$ is prime, $I=A\setminus R$ is an ideal. The hypothesis gives $D_P(Q)\cap I=\emptyset$. By Theorem \ref{theorem_FP}, there exists $Z\in X(A)$ such that $D_P(Q)\subseteq Z$, $(P,Z)\in R_A$ and $Z\cap I=\emptyset$. Since $Q\subseteq D_P(Q)$, we obtain $Q\subseteq Z\subseteq R$.

$(2)\Rightarrow(1)$: By the minimality of $D_P(Q)$, the conditions $Q\subseteq Z$ and $(P,Z)\in R_A$ imply $D_P(Q)\subseteq Z\subseteq R$.
\end{proof}

\begin{lemma} \label{lema_aux2} \cite[Proposition 11(5)]{Castro2}
Let $\langle A, \Rightarrow \rangle \in {\bf{GI}}$. Let $F \in {\rm{Fi}}(A)$ and $Q \in X(A)$ such that $(F,Q) \in R_{A}$. Then there exists $P \in X(A)$ such that $(P,Q) \in R_{A}$ and $F \subseteq P$.
\end{lemma}
\begin{proof}
We take
\[
\Sigma = \{ H \in {\rm{Fi}}(A) \colon (H,Q) \in R_{A} \text{ and } F \subseteq H \}.
\]
Since $F \in \Sigma$, the set $\Sigma$ is non-empty. If $\mathcal C\subseteq\Sigma$ is a non-empty chain, then $\bigcup\mathcal C$ is a filter containing $F$. If $t\in(a\Rightarrow b)\cap\bigcup\mathcal C$ and $a\in Q$, then $t$ belongs to a member $H\in\mathcal C$, and $(H,Q)\in R_A$ yields $b\in Q$. Thus $(\bigcup\mathcal C,Q)\in R_A$, so every chain has an upper bound in $\Sigma$, with $F$ serving for the empty chain. Zorn's lemma therefore gives $P \in {\rm{Fi}}(A)$ such that $P$ is maximal in $\Sigma$.
If $P=A$, then $(1 \Rightarrow 0) \cap P \neq \emptyset$, and since $1 \in Q$ and $(P,Q) \in R_{A}$, we get $0 \in Q$, i.e., $Q=A$, which is a contradiction. So, $P \neq A$.
Now we prove that $P \in X(A)$. Suppose there exist $a,b \in A$ such that $a \vee b \in P$ and $a,b \notin P$. Consider $P_{a}={\rm{Fig}}(P \cup \{a\})$ and $P_{b}={\rm{Fig}}(P \cup \{b\})$. Since $F \subseteq P \subset P_{a}$, $F \subseteq P \subset P_{b}$, and $P$ is maximal in $\Sigma$, we have $(P_{a},Q) \notin R_{A}$ and $(P_{b},Q) \notin R_{A}$. Then there exist $x_{1}, x_{2}, y_{1}, y_{2} \in A$ such that $(x_{1} \Rightarrow y_{1}) \cap P_{a} \neq \emptyset$, $(x_{2} \Rightarrow y_{2}) \cap P_{b} \neq \emptyset$, $x_{1}, x_{2} \in Q$ and $y_{1}, y_{2} \notin Q$. Let $x = x_{1} \wedge x_{2} \in Q$. As $Q$ is a prime filter, $y = y_{1} \vee y_{2} \notin Q$. Since $x \leq x_{1}, x_{2}$, we have $(x \Rightarrow y_{1}) \cap P_{a} \neq \emptyset$ and $(x \Rightarrow y_{2}) \cap P_{b} \neq \emptyset$. Then there is $l_{1} \in P_{a}$ such that $l_{1} \in (x \Rightarrow y_{1}) \subseteq x \Rightarrow y$ and there is $l_{2} \in P_{b}$ such that $l_{2} \in (x \Rightarrow y_{2}) \subseteq x \Rightarrow y$. So, there exist $p_{1}, p_{2} \in P$ such that $p_{1} \wedge a \leq l_{1}$ and $p_{2} \wedge b \leq l_{2}$. If we take $p = p_{1} \wedge p_{2} \in P$, then $p \wedge a, p \wedge b \in x \Rightarrow y$, which implies $(p \wedge a) \vee (p \wedge b) = p \wedge (a \vee b) \in x \Rightarrow y$. So, $(x \Rightarrow y) \cap P \neq \emptyset$, and since $x \in Q$ and $(P,Q) \in R_{A}$, we get $y \in Q$, which is a contradiction. We conclude that $P \in X(A)$.
\end{proof}

\section{Commutative gi-lattices} \label{sec3}

We introduce commutative gi-lattices as an extension of the self-distributive weak Heyting algebras introduced in \cite{Ghorbani2}. We study their connection with other subclasses of gi-lattices and present some characterizations.

\begin{definition}
Let $\langle A, \Rightarrow \rangle \in {\bf{GI}}$. We say that $\langle A, \Rightarrow \rangle$ is commutative if the following identity is satisfied for all $a,b,c \in A$:
\[
a \Rightarrow (b \Rightarrow c) = b \Rightarrow (a \Rightarrow c).
\]
\end{definition}

We denote by ${\bf{CGI}}$ the class of commutative gi-lattices.

\begin{lemma} \label{lema_GIC1}
Let $\langle A, \Rightarrow \rangle \in {\bf{CGI}}$. Let $a,b \in A$. Then:
\begin{enumerate}
\item $a \Rightarrow (b \Rightarrow a) = A$,
\item $a \Rightarrow (b \Rightarrow (a \wedge b)) = A$,
\item $a \Rightarrow b \subseteq 1 \Rightarrow (a \Rightarrow b)$.
\end{enumerate}
\end{lemma}
\begin{proof}
$(1)$ Since $\langle A, \Rightarrow \rangle$ is commutative, by (GI4) of Definition \ref{def: generalized implication} we have
\[
a \Rightarrow (b \Rightarrow a) = b \Rightarrow (a \Rightarrow a) = b \Rightarrow A = A.
\]

$(2)$ By (GI1) and (GI4) of Definition \ref{def: generalized implication} and item (1), it follows that
\[
a \Rightarrow (b \Rightarrow (a \wedge b)) = a \Rightarrow ((b \Rightarrow a) \cap (b \Rightarrow b)) = a \Rightarrow (b \Rightarrow a) = A.
\]

$(3)$ Suppose there exist $a,b \in A$ such that $a \Rightarrow b \nsubseteq 1 \Rightarrow (a \Rightarrow b)$, i.e., there is $x \in a \Rightarrow b$ such that $x \notin 1 \Rightarrow (a \Rightarrow b)$. Then there exists $P \in X(A)$ such that $x \in P$ and $P \cap (1 \Rightarrow (a \Rightarrow b)) = \emptyset$. Since $\langle A, \Rightarrow \rangle$ is commutative, $P \cap (a \Rightarrow (1 \Rightarrow b)) = \emptyset$. So, by Lemma \ref{lema_F to I}(4), there exists $Q \in X(A)$ such that $(P,Q) \in R_{A}$, $a \in Q$ and $Q \cap (1 \Rightarrow b) = \emptyset$. Since $x \in P \cap (a \Rightarrow b)$, $a \in Q$ and $(P,Q) \in R_{A}$, it follows that $b \in Q$. On the other hand, by item (1), $b \Rightarrow (1 \Rightarrow b) = A$, so $P \cap (b \Rightarrow (1 \Rightarrow b)) \neq \emptyset$ and there is $y \in P$ such that $y \in b \Rightarrow (1 \Rightarrow b)$. By (\ref{a to I}), there exists $i \in 1 \Rightarrow b$ such that $y \in b \Rightarrow i$. Then $y \in P \cap (b \Rightarrow i)$, $b \in Q$ and $(P,Q) \in R_{A}$ imply $i \in Q$. Therefore, $i \in Q \cap (1 \Rightarrow b)$, which is a contradiction. So, $a \Rightarrow b \subseteq 1 \Rightarrow (a \Rightarrow b)$.
\end{proof}

\begin{example} \label{example_principal} 
In the principal-ideal construction of Example \ref{Example_2},
\[
a\Rightarrow(b\Rightarrow c)
=\bigcup_{i\leq b\to c}(a\to i]
=(a\to(b\to c)].
\]
Consequently, a weak Heyting algebra is commutative if and only if its associated gi-lattice is commutative. Conversely, if every value $a\Rightarrow b$ is principal, its greatest element defines a unique WH-operation on $A$, and (GI1)--(GI4) become the WH-identities. In particular, every finite gi-lattice arises in this way.
\end{example}

Our next goal is to compare the class of commutative gi-lattices with other subclasses of gi-lattices.

\begin{lemma} \label{lemma_SC1}
 ${\bf{CGI}} \subseteq {\bf{TGI}}$.
\end{lemma}
\begin{proof}
By Lemma \ref{lema_GIC1}(3), $a\Rightarrow b\subseteq1\Rightarrow(a\Rightarrow b)$. Since $c\leq1$, antitonicity in the first coordinate gives $1\Rightarrow(a\Rightarrow b)\subseteq c\Rightarrow(a\Rightarrow b)$ for every $c\in A$. This is the defining condition of $ {\bf TGI}$.
\end{proof}

\begin{remark} (${\bf{CGI}} \subset {\bf{TGI}}$): \label{remark_CGI}
Let $A$ be the three-element chain $0 < a < 1$, and consider the following generalized implication:
\vspace{0.4cm}
\begin{center}
\begin{tabular}{c|ccc}
$\Rightarrow$ & $0$ & $a$ & $1$\tabularnewline
\hline
$0$ & $A$     & $A$          & $A$\tabularnewline
$a$ & $\{0\}$ & $A$          & $A$\tabularnewline
$1$ & $\{0\}$ & $\{0\}$     & $A$\tabularnewline
\end{tabular}
\par
\end{center}
\vspace{0.4cm}
Then $\langle A, \Rightarrow \rangle \in {\bf{TGI}}$ and $\langle A, \Rightarrow \rangle \notin {\bf{CGI}}$. Indeed,
\[
a \Rightarrow (1 \Rightarrow a) = \{0\} \neq A = 1 \Rightarrow (a \Rightarrow a).
\]
\end{remark}

\begin{remark} \label{remark_constant_not_CGI}
The gi-lattice of Remark \ref{remark_CGI} is the instance of Example \ref{Example_constant} with $A$ the three-element chain and $I = \{0\}$, and the computation given there generalizes as follows. Let $\langle A, \Rightarrow_{I} \rangle$ be as in Example \ref{Example_constant}, and suppose that there exists $a \in A$ with $a \neq 1$ and $a \notin I$. Then $1 \Rightarrow_{I} a = I$ and, by (\ref{a to I}),
\[
a \Rightarrow_{I} (1 \Rightarrow_{I} a) = \bigcup \{ a \Rightarrow_{I} i \colon i \in I \} = I,
\]
because $a \nleq i$ for every $i \in I$ (otherwise $a \in I$, since $I$ is a downset). On the other hand, $1 \Rightarrow_{I} (a \Rightarrow_{I} a) = 1 \Rightarrow_{I} A = A$. Hence $\langle A, \Rightarrow_{I} \rangle \notin {\bf{CGI}}$ whenever $A \neq I \cup \{1\}$. Note that $\langle A, \Rightarrow_{I} \rangle \in {\bf{TGI}}$: if $x \leq y$ then $z \Rightarrow_{I} (x \Rightarrow_{I} y) = z \Rightarrow_{I} A = A$, and if $x \nleq y$ then $x \Rightarrow_{I} y = I \subseteq z \Rightarrow_{I} I$, since $0 \in I$ and $z \Rightarrow_{I} 0 \in \{A, I\}$. In particular, the gi-lattice $\langle {\rm{FC}}(\mathbb{N}), \Rightarrow_{{\rm{Fin}}} \rangle$ of Example \ref{Example_constant} is an infinite member of ${\bf{TGI}} \setminus {\bf{CGI}}$ whose implication values are not principal: it suffices to take $a = \mathbb{N} \setminus \{0\}$, which is cofinite and different from $\mathbb{N}$.
\end{remark}

\begin{remark} (${\bf{DL}} \subset {\bf{CGI}}$): \label{remark_DL}
If $\Rightarrow$ is the relative annihilator, then for all $a,b,c \in A$ we have, by (\ref{a to I}), that $x \in a \Rightarrow (b \Rightarrow c)$ if and only if there is $i \in A$ with $i \wedge b \leq c$ and $x \wedge a \leq i$, if and only if $x \wedge a \wedge b \leq c$. This condition is symmetric in $a$ and $b$, so ${\bf{DL}} \subseteq {\bf{CGI}}$. To see that the inclusion is proper, let $A$ be the three-element chain $0 < a < 1$, and consider the following generalized implication:
\vspace{0.4cm}
\begin{center}
\begin{tabular}{c|ccc}
$\Rightarrow$ & $0$ & $a$ & $1$\tabularnewline
\hline
$0$ & $A$     & $A$          & $A$\tabularnewline
$a$ & $A$     & $A$          & $A$\tabularnewline
$1$ & $\{0\}$ & $\{0\}$     & $A$\tabularnewline
\end{tabular}
\par
\end{center}
\vspace{0.4cm}
Then $\langle A, \Rightarrow \rangle \in {\bf{CGI}}$ and $\Rightarrow$ is not the relative annihilator because $\langle a, 0 \rangle = \{0\}$, which is not equal to $a \Rightarrow 0 = A$.
\end{remark}

\begin{remark} (${\bf{CGI}} \neq {\bf{BasGI}}$): \label{remark_CGIBasGI}
Consider the gi-lattice $\langle A, \Rightarrow \rangle$ from Remark \ref{remark_DL}. Then $\langle A, \Rightarrow \rangle \in {\bf{CGI}}$. Since $(a] \nsubseteq 1 \Rightarrow a = \{0\}$, we have $\langle A, \Rightarrow \rangle \notin {\bf{BasGI}}$. On the other hand, consider the three-element chain $A$ given by $0 < a < 1$ and the following generalized implication:
\vspace{0.4cm}
\begin{center}
\begin{tabular}{c|ccc}
$\Rightarrow$ & $0$ & $a$ & $1$\tabularnewline
\hline
$0$ & $A$     & $A$          & $A$\tabularnewline
$a$ & $(a]$   & $A$          & $A$\tabularnewline
$1$ & $(a]$    & $(a]$       & $A$\tabularnewline
\end{tabular}
\par
\end{center}
\vspace{0.4cm}
Then $\langle A, \Rightarrow \rangle \in {\bf{BasGI}}$ and $\langle A, \Rightarrow \rangle \notin {\bf{CGI}}$. Indeed,
\[
1 \Rightarrow (a \Rightarrow 0) = (a] \neq A = a \Rightarrow (1 \Rightarrow 0).
\]
\end{remark}

\begin{remark} (${\bf{CGI}} \neq {\bf{TRGI}}$): \label{remark_CGITRGI}
Consider the gi-lattice $\langle A, \Rightarrow \rangle$ from Remark \ref{remark_CGI}. Then $\langle A, \Rightarrow \rangle \notin {\bf{CGI}}$, and it is easy to see that $\langle A, \Rightarrow \rangle \in {\bf{TRGI}}$. On the other hand, let $A$ be the three-element chain $0 < a < 1$ and consider the following generalized implication:
\vspace{0.4cm}
\begin{center}
\begin{tabular}{c|ccc}
$\Rightarrow$ & $0$ & $a$ & $1$\tabularnewline
\hline
$0$ & $A$     & $A$          & $A$\tabularnewline
$a$ & $A$   & $A$          & $A$\tabularnewline
$1$ & $(a]$    & $(a]$       & $A$\tabularnewline
\end{tabular}
\par
\end{center}
\vspace{0.4cm}
It follows that $\langle A, \Rightarrow \rangle \in {\bf{CGI}}$ and $\langle A, \Rightarrow \rangle \notin {\bf{TRGI}}$, since $(a] \cap (a \Rightarrow 0) = (a] \nsubseteq \{0\}$.
\end{remark}

\begin{remark} \label{remark_CGIcapBas}(${\bf{CGI}} \cap {\bf{BasGI}} \neq {\bf{DL}}$):
Let $A$ be the three-element chain $0 < a < 1$ and consider the following generalized implication:
\vspace{0.4cm}
\begin{center}
\begin{tabular}{c|ccc}
$\Rightarrow$ & $0$ & $a$ & $1$\tabularnewline
\hline
$0$ & $A$         & $A$          & $A$\tabularnewline
$a$ & $\{0\}$     & $A$          & $A$\tabularnewline
$1$ & $\{0\}$    & $A$           & $A$\tabularnewline
\end{tabular}
\par
\end{center}
\vspace{0.4cm}
Then $\langle A, \Rightarrow \rangle \in {\bf{CGI}} \cap {\bf{BasGI}}$. Furthermore, $\Rightarrow$ is not the relative annihilator because $\langle 1, a \rangle = (a]$, which does not equal $1 \Rightarrow a = A$.
\end{remark}

\begin{lemma} \label{lemma_SC2}
Let $\langle A, \Rightarrow \rangle \in {\bf{GI}}$. Then the generalized implication $\Rightarrow$ is the relative annihilator if and only if $\langle A, \Rightarrow \rangle \in {\bf{CGI}} \cap {\bf{RGI}}$.
\end{lemma}
\begin{proof}
If $\Rightarrow$ is the relative annihilator, membership in $ {\bf CGI}\cap{\bf RGI}$ follows from Remarks \ref{remark_DL} and \ref{remark_DL0}.

Conversely, let $\langle A,\Rightarrow\rangle\in{\bf CGI}\cap{\bf RGI}$ and $a\in A$. By commutativity,
\[
a\Rightarrow(1\Rightarrow a)=1\Rightarrow(a\Rightarrow a)=A.
\]
Since $1\in a\Rightarrow(1\Rightarrow a)$, equation (\ref{a to I}) gives an $i\in1\Rightarrow a$ with $1\in a\Rightarrow i$. The $ {\bf RGI}$ condition yields $a=1\wedge a\leq i$. As $1\Rightarrow a$ is an ideal, $a\in1\Rightarrow a$, so $(a]\subseteq1\Rightarrow a$. Conversely, $ {\bf RGI}$ gives $1\Rightarrow a\subseteq(a]$. Thus $1\Rightarrow a=(a]$ for every $a$, and Proposition \ref{propo_annihilator} applies.
\end{proof}

The finite examples given in Remarks \ref{remark_DL}, \ref{remark_CGITRGI} and \ref{remark_CGIcapBas} are instances of the following general construction, which also produces infinite commutative gi-lattices whose implication values are not principal.

\begin{proposition} \label{prop_PJ}
Let $A$ be a bounded distributive lattice, $P \in X(A)$ and $J \in {\rm{Id}}(A)$. Define $\Rightarrow_{P,J} \colon A \times A \to {\rm{Id}}(A)$ by
\[
a \Rightarrow_{P,J} b =
\begin{cases}
J & \text{if } a \in P \text{ and } b \notin P,\\
A & \text{otherwise.}
\end{cases}
\]
Then $\langle A, \Rightarrow_{P,J} \rangle \in {\bf{CGI}}$. Moreover, if $J \cap P = \emptyset$, then $a \Rightarrow_{P,J} (b \Rightarrow_{P,J} c) = (a \wedge b) \Rightarrow_{P,J} c$ for all $a,b,c \in A$, and if $J \cap P \neq \emptyset$, then $a \Rightarrow_{P,J} (b \Rightarrow_{P,J} c) = A$ for all $a,b,c \in A$.
\end{proposition}
\begin{proof}
Throughout the proof we write $\Rightarrow$ instead of $\Rightarrow_{P,J}$, and for $a,b \in A$ we write $\delta(a,b)$ to mean that $a \in P$ and $b \notin P$, so that $a \Rightarrow b = J$ if $\delta(a,b)$ holds and $a \Rightarrow b = A$ otherwise. Since $P$ is a prime filter, for all $a,b,c \in A$ we have:
\begin{enumerate}
\item[(i)] $\delta(a, b \wedge c)$ holds if and only if $\delta(a,b)$ or $\delta(a,c)$ holds, because $b \wedge c \notin P$ if and only if $b \notin P$ or $c \notin P$;
\item[(ii)] $\delta(a \vee b, c)$ holds if and only if $\delta(a,c)$ or $\delta(b,c)$ holds, because $a \vee b \in P$ if and only if $a \in P$ or $b \in P$;
\item[(iii)] if $\delta(a,c)$ holds, then $\delta(a,b)$ or $\delta(b,c)$ holds, according to whether $b \notin P$ or $b \in P$;
\item[(iv)] $\delta(a,a)$ does not hold.
\end{enumerate}
We also note that the intersection of two members of $\{A, J\}$ is $J$ if and only if at least one of them is $J$. Now (GI1) follows from (i): $a \Rightarrow (b \wedge c) = J$ if and only if $\delta(a,b)$ or $\delta(a,c)$ holds, if and only if $(a \Rightarrow b) \cap (a \Rightarrow c) = J$; otherwise both sides are equal to $A$. Condition (GI2) follows from (ii) in the same way. For (GI3), if $a \Rightarrow c = A$ there is nothing to prove; if $a \Rightarrow c = J$, then $\delta(a,c)$ holds, so by (iii) at least one of $a \Rightarrow b$, $b \Rightarrow c$ is equal to $J$, and $(a \Rightarrow b) \cap (b \Rightarrow c) \subseteq J$. Finally, (GI4) is (iv). Hence $\langle A, \Rightarrow \rangle \in {\bf{GI}}$.

Let $I \in {\rm{Id}}(A)$ and $a \in A$. By (\ref{a to I}), $a \Rightarrow I = \bigcup \{ a \Rightarrow i \colon i \in I \}$. If $a \notin P$, then $a \Rightarrow i = A$ for every $i \in I$, so $a \Rightarrow I = A$. If $a \in P$, then $a \Rightarrow i = A$ when $i \in P$ and $a \Rightarrow i = J$ when $i \notin P$; since $0 \in I$, we obtain $a \Rightarrow I = A$ if $I \cap P \neq \emptyset$, and $a \Rightarrow I = J$ if $I \cap P = \emptyset$. As $1 \in P$, in particular $a \Rightarrow A = A$ for every $a \in A$.

Let $a,b,c \in A$. If $\delta(b,c)$ does not hold, then $b \Rightarrow c = A$ and $a \Rightarrow (b \Rightarrow c) = A$. If $\delta(b,c)$ holds, then $b \Rightarrow c = J$, and by the previous paragraph $a \Rightarrow (b \Rightarrow c) = J$ when $a \in P$ and $J \cap P = \emptyset$, while $a \Rightarrow (b \Rightarrow c) = A$ otherwise. Therefore, if $J \cap P \neq \emptyset$, then $a \Rightarrow (b \Rightarrow c) = A$ for all $a,b,c \in A$; and if $J \cap P = \emptyset$, then $a \Rightarrow (b \Rightarrow c) = J$ if and only if $a \in P$, $b \in P$ and $c \notin P$, i.e., if and only if $\delta(a \wedge b, c)$ holds, so $a \Rightarrow (b \Rightarrow c) = (a \wedge b) \Rightarrow c$. In both cases $a \Rightarrow (b \Rightarrow c)$ is symmetric in $a$ and $b$, and thus $\langle A, \Rightarrow \rangle \in {\bf{CGI}}$.
\end{proof}

\begin{remark} \label{remark_PJ_instances}
For the three-element chain $0 < a < 1$, the gi-lattices of Remarks \ref{remark_DL}, \ref{remark_CGITRGI} and \ref{remark_CGIcapBas} are $\Rightarrow_{P,J}$ with $(P,J) = (\{1\}, \{0\})$, $(\{1\}, (a])$ and $([a), \{0\})$, respectively, and the gi-lattice of Example \ref{gi_sd} below is $\Rightarrow_{P,J}$ with $P = [b)$ and $J = (b]$. The commutative gi-lattice of Example \ref{example_theta} is not of this form, since its non-trivial values $\{0\}$ and $(a]$ are different.
\end{remark}

\begin{example} \label{example_infinite}
Let $A = {\rm{FC}}(\mathbb{N})$ be the Boolean algebra of all finite and cofinite subsets of $\mathbb{N}$, let $P$ be the prime filter of all cofinite subsets of $\mathbb{N}$, and let $J = {\rm{Fin}} = A \setminus P$ be the ideal of all finite subsets of $\mathbb{N}$. By Proposition \ref{prop_PJ}, the map
\[
X \Rightarrow Y =
\begin{cases}
{\rm{Fin}} & \text{if } X \text{ is cofinite and } Y \text{ is finite},\\
A & \text{otherwise}
\end{cases}
\]
is a generalized implication on $A$ and $\langle A, \Rightarrow \rangle \in {\bf{CGI}}$. Since $J \cap P = \emptyset$, it also satisfies $X \Rightarrow (Y \Rightarrow Z) = (X \cap Y) \Rightarrow Z$ for all $X,Y,Z \in A$. The ideal ${\rm{Fin}}$ has no greatest element, because if $F$ is a finite subset of $\mathbb{N}$ and $n \notin F$, then $F \cup \{n\}$ is a finite subset strictly above $F$. Thus $\mathbb{N} \Rightarrow \emptyset = {\rm{Fin}}$ is not a principal ideal, and by Example \ref{example_principal} the gi-lattice $\langle A, \Rightarrow \rangle$ does not arise from any WH-algebra. Moreover:
\begin{enumerate}
\item $\langle A, \Rightarrow \rangle \in {\bf{BasGI}}$: if $X$ is cofinite, then $1 \Rightarrow X = \mathbb{N} \Rightarrow X = A$, and if $X$ is finite, then $1 \Rightarrow X = {\rm{Fin}} \supseteq (X]$.
\item $\langle A, \Rightarrow \rangle \notin {\bf{RGI}}$, since $(\mathbb{N}] \cap (\mathbb{N} \Rightarrow \emptyset) = {\rm{Fin}} \nsubseteq (\emptyset] = \{\emptyset\}$. Hence, by Lemma \ref{lemma_SC2}, $\Rightarrow$ is not the relative annihilator; this can also be seen directly, since $\langle \mathbb{N}, \emptyset \rangle = \{\emptyset\} \neq \mathbb{N} \Rightarrow \emptyset$.
\end{enumerate}
Therefore $\langle A, \Rightarrow \rangle$ is an infinite member of $({\bf{CGI}} \cap {\bf{BasGI}}) \setminus {\bf{DL}}$ whose implication values are not principal, which is the infinite counterpart of Remark \ref{remark_CGIcapBas}. The same construction applies to any Boolean algebra $B$ with a non-principal ultrafilter $P$, taking $J = B \setminus P$, and more generally to any bounded distributive lattice having a prime filter whose complement is a non-principal ideal.
\end{example}

\begin{remark}
Based on Remarks \ref{remark_CGI}, \ref{remark_DL}, \ref{remark_CGIBasGI}, \ref{remark_CGITRGI} and Lemmas \ref{lemma_SC1} and \ref{lemma_SC2}, we obtain the following diagram of inclusions among subclasses of gi-lattices:
\vspace{0.4cm}
\begin{center}
\begin{tikzpicture}
\filldraw [thick]
(0,1) circle (2pt) node[left] {${\bf{BasGI}}$} --
(1,0) circle (2pt) node[right] {${\bf{DL}}$}--(1,1) circle (2pt) node[left] {${\bf{CGI}}$}--(1,2)
(1,0)--(2,1) circle (2pt) node[right] {${\bf{TRGI}}$} --
(1,2) circle (2pt) node[left] {${\bf{TGI}}$}--(0,1)
(2,1)--(3,2)circle (2pt) node[right] {${\bf{RGI}}$}--
(2,3) circle (2pt) node[left]{${\bf{GI}}$} -- (1,2);
\end{tikzpicture}
\end{center}
\vspace{0.2cm}
Note that, by Lemma \ref{lemma_SC2}, ${\bf{CGI}} \cap {\bf{RGI}} = {\bf{CGI}} \cap {\bf{TRGI}} = {\bf{DL}}$, whereas ${\bf{CGI}} \cap {\bf{BasGI}} \neq {\bf{DL}}$ by Remark \ref{remark_CGIcapBas}. Hence the diagram only displays inclusions and should not be read as a lattice in which meets are intersections. By Remark \ref{remark_constant_not_CGI} and Example \ref{example_infinite}, the inclusion ${\bf{CGI}} \subset {\bf{TGI}}$ and the inequality ${\bf{CGI}} \cap {\bf{BasGI}} \neq {\bf{DL}}$ are also witnessed by infinite gi-lattices whose implication values are not principal.
\end{remark}

The following characterizations extend the prime-filter characterizations for commutative weak Heyting algebras in \cite{Celani-Hernan}. The extension replaces element-valued implication by ideal-valued implication; it retains the same relational pattern.

\begin{theorem} \label{theorem_p1}
Let $\langle A, \Rightarrow \rangle \in {\bf{GI}}$. Then the following conditions are equivalent:
\begin{enumerate}
\item $\langle A, \Rightarrow \rangle \in {\bf{CGI}}$.
\item For all $P,Q,R \in X(A)$, if $(P,Q) \in R_{A}$ and $(Q,R) \in R_{A}$, then $D_{D_{P}(R)}(Q) \subseteq R$.
\end{enumerate}
\end{theorem}
\begin{proof}
$(1)$ implies $(2)$: Let $P,Q,R \in X(A)$ be such that $(P,Q) \in R_{A}$ and $(Q,R) \in R_{A}$. Let $a \in D_{D_{P}(R)}(Q)$. Then there is $q \in Q$ such that $D_{P}(R) \cap (q \Rightarrow a) \neq \emptyset$. So, there is $x \in A$ such that $x \in D_{P}(R) \cap (q \Rightarrow a)$, i.e., there is $r \in R$ such that $P \cap (r \Rightarrow x) \neq \emptyset$. Hence there exists $y \in P$ such that $y \in r \Rightarrow x$. Since $r \in [r)$ and $x \in q \Rightarrow a$, it follows that $y \in r \Rightarrow (q \Rightarrow a)$. Thus, $P \cap (r \Rightarrow (q \Rightarrow a)) \neq \emptyset$. Since $\langle A, \Rightarrow \rangle$ is commutative, $P \cap (q \Rightarrow (r \Rightarrow a)) \neq \emptyset$. Then there is $z \in P$ such that $z \in q \Rightarrow (r \Rightarrow a)$, and by Lemma \ref{lema_F to I}, there exist $f \in [q)$ and $i \in r \Rightarrow a$ such that $z \in f \Rightarrow i$. Then $f \in Q$. Since $P \cap (f \Rightarrow i) \neq \emptyset$, $f \in Q$ and $(P,Q) \in R_{A}$\bl{,} we have $i \in Q$. It follows that $i \in Q \cap (r \Rightarrow a)$, i.e., $Q \cap (r \Rightarrow a) \neq \emptyset$. Since $r \in R$ and $(Q,R) \in R_{A}$, we get $a \in R$. Therefore, $D_{D_{P}(R)}(Q) \subseteq R$.

$(2)$ implies $(1)$: Suppose there exist $a,b,c \in A$ such that $a \Rightarrow (b \Rightarrow c) \nsubseteq b \Rightarrow (a \Rightarrow c)$. So, there is $x \in a \Rightarrow (b \Rightarrow c)$ such that $x \notin b \Rightarrow (a \Rightarrow c)$. Then there is $P \in X(A)$ such that $x \in P$ and $P \cap (b \Rightarrow (a \Rightarrow c)) = \emptyset$. By Lemma \ref{lema_F to I}, there exists $Q \in X(A)$ such that $(P,Q) \in R_{A}$, $b \in Q$ and $Q \cap (a \Rightarrow c) =  \emptyset$. Thus, by Corollary \ref{corollary_FP}, there is $R \in X(A)$ such that $(Q,R) \in R_{A}$, $a \in R$ and $c \notin R$. Since $(P,Q) \in R_{A}$ and $(Q,R) \in R_{A}$, by hypothesis we have $D_{D_{P}(R)}(Q) \subseteq R$. On the other hand, we show that $c \in D_{D_{P}(R)}(Q)$. Indeed, since $a \in R$, we have $x \in R \Rightarrow (b \Rightarrow c)$ by Lemma \ref{lema_F to I}(1), so $P \cap (R \Rightarrow (b \Rightarrow c)) \neq \emptyset$, and by Lemma \ref{lema_F to I}(3), we have $D_{P}(R) \cap (b \Rightarrow c) \neq \emptyset$. By definition and since $b \in Q$, it follows that $c \in D_{D_{P}(R)}(Q)$. Therefore, $c \in R$, which is a contradiction. We conclude that $a \Rightarrow (b \Rightarrow c) \subseteq b \Rightarrow (a \Rightarrow c)$, and interchanging the roles of $a$ and $b$ we obtain the equality.
\end{proof}

\begin{theorem} \label{charact_2}
Let $\langle A, \Rightarrow \rangle \in {\bf{GI}}$. Then the following conditions are equivalent:
\begin{enumerate}
\item $\langle A, \Rightarrow \rangle \in {\bf{CGI}}$.
\item For all $P,Q,R \in X(A)$, if $(P,Q) \in R_{A}$ and $(Q,R) \in R_{A}$, then there exist $F,K,H \in X(A)$ such that $(P,F) \in R_{A}$, $(H,K) \in R_{A}$ and $Q \subseteq K \subseteq R \subseteq F \subseteq H$.
\end{enumerate}
\end{theorem}
\begin{proof}
$(1)$ implies $(2)$: Let $P,Q,R \in X(A)$ be such that $(P,Q) \in R_{A}$ and $(Q,R) \in R_{A}$. By Theorem \ref{theorem_p1}, we have $D_{D_{P}(R)} (Q) \subseteq R$. So, $D_{D_{P}(R)} (Q) \cap R^{c} =  \emptyset$, where $R^{c} = A \setminus R$ is an ideal because $R$ is a prime filter. By Theorem \ref{theorem_FP}, there exists $K \in X(A)$ such that $(D_{P}(R), K) \in R_{A}$ and $Q \subseteq K \subseteq R$. Then, by Lemma \ref{lema_aux2}, there is $H \in X(A)$ such that $(H,K) \in R_{A}$ and $D_{P}(R) \subseteq H$. By Lemma \ref{lema_aux1}, there is $F \in X(A)$ such that $R \subseteq F \subseteq H$ and $(P,F) \in R_{A}$. In summary, there exist $F,K,H \in X(A)$ such that $(P,F) \in R_{A}$, $(H,K) \in R_{A}$ and $Q \subseteq K \subseteq R \subseteq F \subseteq H$.

$(2)$ implies $(1)$: Suppose there exist $a,b,c \in A$ such that $a \Rightarrow (b \Rightarrow c) \nsubseteq b \Rightarrow (a \Rightarrow c)$. So, there is $x \in a \Rightarrow (b \Rightarrow c)$ such that $x \notin b \Rightarrow (a \Rightarrow c)$. Then there is $P \in X(A)$ such that $x \in P$ and $(b \Rightarrow (a \Rightarrow c)) \cap P = \emptyset$. By Lemma \ref{lema_F to I}, there is $Q \in X(A)$ such that $(P,Q) \in R_{A}$, $b \in Q$ and $(a \Rightarrow c) \cap Q = \emptyset$. Applying Lemma \ref{lema_F to I} once again, there is $R \in X(A)$ such that $(Q,R) \in R_{A}$, $a \in R$ and $c \notin R$. So, by hypothesis, there exist $F,K,H \in X(A)$ such that $(P,F) \in R_{A}$, $(H,K) \in R_{A}$ and $Q \subseteq K \subseteq R \subseteq F \subseteq H$. Then $a \in R \subseteq F$, i.e., $a \in F$. On the other hand, since $x \in a \Rightarrow (b \Rightarrow c)$, by (\ref{a to I}) there is $i \in b \Rightarrow c$ such that $x \in a \Rightarrow i$. Thus, $x \in (a \Rightarrow i) \cap P$, i.e., $(a \Rightarrow i) \cap P \neq \emptyset$. Since $a \in F$ and $(P,F) \in R_{A}$, we get $i \in F$. So, $i \in (b \Rightarrow c) \cap F$, i.e., $(b \Rightarrow c) \cap F \neq \emptyset$. As $F \subseteq H$, it follows that $(b \Rightarrow c) \cap H \neq \emptyset$. Moreover, $b \in Q \subseteq K$, and since $(H,K) \in R_{A}$, we get $c \in K \subseteq R$. We conclude that $c \in R$, which is a contradiction. Then $a \Rightarrow (b \Rightarrow c) \subseteq b \Rightarrow (a \Rightarrow c)$. The other inclusion is similar and therefore $\langle A, \Rightarrow \rangle \in {\bf{CGI}}$.
\end{proof}

\begin{corollary} \label{coro_aa}
Let $\langle A, \Rightarrow \rangle \in {\bf{CGI}}$ and $a,b \in A$. Then $a \Rightarrow (a \Rightarrow b) = 1 \Rightarrow (a \Rightarrow b)$.
\end{corollary}
\begin{proof}
Since $a\leq1$, antitonicity gives $1\Rightarrow(a\Rightarrow b)\subseteq a\Rightarrow(a\Rightarrow b)$. For the converse, suppose that $x\in a\Rightarrow(a\Rightarrow b)$ but $x\notin1\Rightarrow(a\Rightarrow b)$. By prime-filter separation there exists $P\in X(A)$ with $x\in P$ and
\[
P\cap(1\Rightarrow(a\Rightarrow b))=\emptyset.
\]
By commutativity, $P\cap(a\Rightarrow(1\Rightarrow b))=\emptyset$. Lemma \ref{lema_F to I}(4) therefore gives $Q\in X(A)$ such that
\[
(P,Q)\in R_A,\qquad a\in Q,\qquad Q\cap(1\Rightarrow b)=\emptyset.
\]
Corollary \ref{corollary_FP}, applied to $1\Rightarrow b$, gives $R\in X(A)$ with $(Q,R)\in R_A$ and $b\notin R$. By Theorem \ref{charact_2}, there exist $F,K,H\in X(A)$ such that
\[
(P,F)\in R_A,\qquad(H,K)\in R_A,\qquad
Q\subseteq K\subseteq R\subseteq F\subseteq H.
\]
On the other hand, $x\in a\Rightarrow(a\Rightarrow b)$ yields an $i\in a\Rightarrow b$ with $x\in a\Rightarrow i$. Since $x\in P$, $a\in Q$ and $(P,Q)\in R_A$, we have $i\in Q\subseteq H$. Thus $H\cap(a\Rightarrow b)\neq\emptyset$. As $a\in Q\subseteq K$ and $(H,K)\in R_A$, it follows that $b\in K\subseteq R$, a contradiction. Hence the reverse inclusion also holds.
\end{proof}

\section{Gi-deductive systems and gi-congruences} \label{sec4}

In \cite{Castro2}, the authors introduced the notion of gi-filter. If $\langle A, \Rightarrow \rangle$ is a gi-lattice, recall that a set $F \subseteq A$ is a gi-filter if it is a lattice filter of $A$ and for every $a \in F$, we have $F \cap (1 \Rightarrow a) \neq \emptyset$. We extend the deductive systems studied for self-distributive weak Heyting algebras in \cite{Ghorbani1} to the ideal-valued setting.

\begin{definition}
Let $\langle A, \Rightarrow \rangle \in {\bf{GI}}$. A set $D \subseteq A$ is a gi-deductive system of $A$ if it satisfies the following conditions for all $a,b \in A$:
\begin{enumerate}
\item[(GDS1)] $1 \in D$,
\item[(GDS2)] $a \in D$ and $D \cap (a \Rightarrow b) \neq \emptyset$ imply $b \in D$.
\end{enumerate}
\end{definition}

\begin{remark} \label{remark_gds_not_gif}
Consider the gi-lattice $\langle A, \Rightarrow \rangle$ of Remark \ref{remark_CGI}. We know that $\langle A, \Rightarrow \rangle \notin {\bf{CGI}}$. If $D = [a)$, then $D$ is a gi-deductive system which is not a gi-filter, since $a \in D$ but $D \cap (1 \Rightarrow a) = \emptyset$. Together with Example \ref{gi_sd} below, this shows that, in the class of gi-lattices, the notions of gi-deductive system and gi-filter are not comparable.
\end{remark}

\begin{lemma} \label{gids implica gif}
Let $\langle A, \Rightarrow \rangle \in {\bf{CGI}}$. Then every gi-deductive system is a gi-filter.
\end{lemma}
\begin{proof}
Let $D$ be a gi-deductive system of $A$. Let $a,b \in A$ be such that $a \in D$ and $a \leq b$. Then $a \Rightarrow b = A$. Since $1 \in D \cap A$, we have $D \cap (a \Rightarrow b) \neq \emptyset$. By (GDS2), we get $b \in D$, so $D$ is an upset. Let $a,b \in D$. By Lemma \ref{lema_GIC1}(2), $a \Rightarrow (b \Rightarrow (a \wedge b)) = A$. Then $1 \in D \cap A$ implies
\[
D \cap (a \Rightarrow (b \Rightarrow (a \wedge b))) \neq \emptyset,
\]
i.e., there is $x \in D$ such that $x \in a \Rightarrow (b \Rightarrow (a \wedge b))$. By (\ref{a to I}), there exists $i \in b \Rightarrow (a \wedge b)$ such that $x \in a \Rightarrow i$. Then $D \cap (a \Rightarrow i) \neq \emptyset$, and by (GDS2), we have $i \in D$. It follows that $D \cap (b \Rightarrow (a \wedge b)) \neq \emptyset$, and again by (GDS2), we have $a \wedge b \in D$. So, $D$ is a filter of $A$.

It remains to show that $D$ is a gi-filter. Let $a \in D$. By Lemma \ref{lema_GIC1}(1) with $b=1$, we have $a \Rightarrow (1 \Rightarrow a) = A$. So, $D \cap (a \Rightarrow (1 \Rightarrow a)) \neq \emptyset$ and there is $x \in D$ such that $x \in a \Rightarrow (1 \Rightarrow a)$. By (\ref{a to I}), there exists $i \in 1 \Rightarrow a$ such that $x \in a \Rightarrow i$. Hence, by (GDS2), $D \cap (a \Rightarrow i) \neq \emptyset$ implies $i \in D$, i.e., $D \cap (1 \Rightarrow a) \neq \emptyset$, and therefore $D$ is a gi-filter.
\end{proof}

\begin{example} \label{gi_sd}
Consider the four-element Boolean lattice $A=\{0,a,b,1\}$, where $a$ and $b$ are incomparable. The following table is the principal-ideal version of \cite[Example 3.11]{Ghorbani1}:
\vspace{0.4cm}
\begin{center}
\begin{tabular}{c|cccc}
$\Rightarrow$ & $0$ & $a$ & $b$ & $1$ \tabularnewline
\hline
$0$ & $A$      & $A$          & $A$   & $A$ \tabularnewline
$a$ & $A$      & $A$          & $A$   & $A$\tabularnewline
$b$ & $(b]$    & $(b]$         & $A$  & $A$\tabularnewline
$1$ & $(b]$    & $(b]$         & $A$  & $A$\tabularnewline
\end{tabular}
\par
\end{center}
\vspace{0.4cm}
Then $\langle A, \Rightarrow \rangle \in {\bf{CGI}}$. It is easy to see that $[b)$ is a gi-filter which is not a gi-deductive system (indeed, $b \in [b)$ and $b \in [b) \cap (b \Rightarrow a)$, but $a \notin [b)$). Note also that $[a)$ is a filter which is not a gi-filter. Indeed, $a \in [a)$ but $[a) \cap (1 \Rightarrow a) = \emptyset$.
\end{example}

\begin{remark} \label{remark_1}
Let $\langle A, \Rightarrow \rangle \in {\bf{GI}}$ and $a,b,c \in A$. If $a \leq b \leq c$, then $c \Rightarrow a = (c \Rightarrow b) \cap (b \Rightarrow a)$. Indeed, the inclusion $\supseteq$ is (GI3), and $\subseteq$ follows from Proposition \ref{prop: properties}.
\end{remark}

The following is an ideal-valued extension of the sufficient direction of \cite[Theorem 3.12]{Ghorbani1}.

\begin{lemma} \label{lemma_gif_gids}
Let $\langle A, \Rightarrow \rangle \in {\bf{CGI}}$. If $1 \Rightarrow a \subseteq (a]$, for all $a \in A$, then the notions of gi-filter and gi-deductive system coincide.
\end{lemma}
\begin{proof}
By Lemma \ref{gids implica gif} we know that every gi-deductive system is a gi-filter. Let $F$ be a gi-filter of $A$, and let $a,b \in A$ such that $a \in F$ and $F \cap (a \Rightarrow b) \neq \emptyset$. Since $b \leq a \vee b \leq 1$, by Remark \ref{remark_1}, (GI2) and (GI4),
\[
1 \Rightarrow b = (1 \Rightarrow (a \vee b)) \cap ((a \vee b) \Rightarrow b) = (1 \Rightarrow (a \vee b)) \cap (a \Rightarrow b).
\]
As $F$ is a gi-filter and $a \in F$, we have $F \cap (1 \Rightarrow a) \neq \emptyset$. So, $F \cap (1 \Rightarrow (a \vee b)) \neq \emptyset$ and there exists $f_{1} \in F$ such that $f_{1} \in 1 \Rightarrow (a \vee b)$. On the other hand, since $F \cap (a \Rightarrow b) \neq \emptyset$, there exists $f_{2} \in F$ such that $f_{2} \in a \Rightarrow b$. Take $f=f_{1} \wedge f_{2} \in F$. Then $f \in (1 \Rightarrow (a \vee b)) \cap (a \Rightarrow b) = 1 \Rightarrow b$ and\bl{,} by hypothesis, $f \in (b]$, i.e., $f \leq b$ and $b \in F$. Hence, $F$ is a gi-deductive system.
\end{proof}

\begin{lemma} \label{gids=gif}
Let $\langle A, \Rightarrow \rangle \in {\bf{GI}}$. Then $\langle A, \Rightarrow \rangle \in {\bf{BasGI}}$ if and only if every filter of $A$ is a gi-filter.
\end{lemma}
\begin{proof}
Assume that $\langle A, \Rightarrow \rangle \in {\bf{BasGI}}$. Let $F$ be a filter of $A$ and $a \in F$. Since $(a] \subseteq 1 \Rightarrow a$, it follows that $F \cap (1 \Rightarrow a) \neq \emptyset$ and $F$ is a gi-filter. Conversely, let $a \in A$. Then, by hypothesis, $[a)$ is a gi-filter and $[a) \cap (1 \Rightarrow a) \neq \emptyset$. So, there exists $f \in [a)$ such that $f \in 1 \Rightarrow a$. If $x \in (a]$, then $x \leq a \leq f$ and $x \in 1 \Rightarrow a$. Hence, $(a] \subseteq 1 \Rightarrow a$ and $\langle A, \Rightarrow \rangle \in {\bf{BasGI}}$.
\end{proof}

\begin{theorem}
Let $\langle A, \Rightarrow \rangle \in {\bf{CGI}}$. Then the generalized implication $\Rightarrow$ is the relative annihilator if and only if every filter is a gi-deductive system.
\end{theorem}
\begin{proof}
Assume that $\Rightarrow$ is the relative annihilator, and let $F$ be a filter of $A$. Then $1 \in F$. Let $a,b \in A$ be such that $a \in F$ and $F \cap (a \Rightarrow b) \neq \emptyset$. Then there exists $x \in F$ such that $x \in a \Rightarrow b = \langle a, b \rangle$, i.e., $x \wedge a \leq b$. Since $a,x \in F$ and $F$ is a filter, we have $a \wedge x \in F$ and $b \in F$. So, $F$ is a gi-deductive system.

Conversely, we prove that $\langle A, \Rightarrow \rangle \in {\bf{RGI}}$. Let $a,b \in A$, and consider $x \in (a] \cap (a \Rightarrow b)$. Then $x \leq a$ and $x \in a \Rightarrow b$. So, $a \Rightarrow b \subseteq x \Rightarrow b$ and $x \in x \Rightarrow b$. Suppose $x \notin (b]$. Then there exists $P \in X(A)$ such that $x \in P$ and $b \notin P$. In particular, $P$ is a filter and\bl{,} by hypothesis, $P$ is a gi-deductive system. Since $x \in P \cap (x \Rightarrow b)$, i.e., $P \cap (x \Rightarrow b) \neq \emptyset$, by (GDS2) we have $b \in P$, which is a contradiction. We conclude that $(a] \cap (a \Rightarrow b) \subseteq (b]$ and $\langle A, \Rightarrow \rangle \in {\bf{CGI}} \cap {\bf{RGI}}$. Then, by Lemma \ref{lemma_SC2}, the generalized implication $\Rightarrow$ is the relative annihilator.
\end{proof}

Recall the notion of gi-congruence introduced in \cite{Castro2}. Let $\langle A, \Rightarrow \rangle \in {\bf{GI}}$ and let $\theta$ be a lattice congruence of $A$. We say that $\theta$ is a gi-congruence if for all $(a,b), (c,d) \in \theta$, the following condition holds:
\[
\forall x \in (a \Rightarrow c) \exists y \in (b \Rightarrow d): (x,y) \in \theta.
\]
For every gi-congruence $\theta$ of $A$, its class $[1]_{\theta}$ is a gi-filter \cite{Castro2}. We study the implication-based relation associated with a gi-filter and its connection with gi-deductive systems. For the principal-ideal case, see \cite{Ghorbani1}; for general results on WH-congruences, see \cite{SanMartin}. The first theorem below holds for all gi-lattices, without assuming commutativity. If $F$ is a filter of $A$, we define the binary relation $\theta(F)$ as follows:
\[
\theta(F) := \{ (a,b) \in A \times A \colon \exists f \in F (f \in (a \Rightarrow b) \cap (b \Rightarrow a)) \}.
\]
Note that $\theta(F) = \{ (a,b) \in A \times A \colon F \cap (a \Rightarrow b) \neq \emptyset \text{ and } F \cap (b \Rightarrow a) \neq \emptyset \}$.

\begin{theorem} \label{theo_congruence}
Let $\langle A, \Rightarrow \rangle \in {\bf{GI}}$ and let $F$ be a gi-filter of $A$. Then $\theta(F)$ is a gi-congruence such that $F \subseteq [1]_{\theta(F)}$. Moreover, if $\psi$ is a gi-congruence of $A$ such that $[1]_{\psi} \subseteq F$, then $\psi \subseteq \theta(F)$.
\end{theorem}
\begin{proof}
It is easy to prove that $\theta(F)$ is reflexive and symmetric. We show that $\theta(F)$ is transitive. Let $(a,b), (b,c) \in \theta(F)$. Then there exist $f_{1}, f_{2} \in F$ such that $f_{1} \in (a \Rightarrow b) \cap (b \Rightarrow a)$ and $f_{2} \in (b \Rightarrow c) \cap (c \Rightarrow b)$. Take $f=f_{1} \wedge f_{2} \in F$. By (GI3) of Definition \ref{def: generalized implication} we have
\[
f \in (a \Rightarrow b) \cap (b \Rightarrow c) \subseteq a \Rightarrow c.
\]
Similarly, $f \in c \Rightarrow a$. So, $f \in (a \Rightarrow c) \cap (c \Rightarrow a)$ and $(a,c) \in \theta(F)$. Thus $\theta(F)$ is an equivalence relation on $A$.

We show that $\theta(F)$ is compatible with $\vee$. If $(a,b), (c,d) \in \theta(F)$, then there exist $f_{1}, f_{2} \in F$ such that $f_{1} \in (a \Rightarrow b) \cap (b \Rightarrow a)$ and $f_{2} \in (c \Rightarrow d) \cap (d \Rightarrow c)$. Since $b \leq b \vee d$, we have $a \Rightarrow b \subseteq a \Rightarrow (b \vee d)$. Similarly, $c \Rightarrow d \subseteq c \Rightarrow (b \vee d)$. If we take $f=f_{1} \wedge f_{2} \in F$, by (GI2) of Definition \ref{def: generalized implication},
\[
f \in (a \Rightarrow b) \cap (c \Rightarrow d) \subseteq (a \Rightarrow (b \vee d)) \cap (c \Rightarrow (b \vee d)) = (a \vee c) \Rightarrow (b \vee d).
\]
Analogously, it is easy to see that $f \in (b \vee d) \Rightarrow (a \vee c)$. Hence, $f \in ((a \vee c) \Rightarrow (b \vee d)) \cap ((b \vee d) \Rightarrow (a \vee c))$ and $(a \vee c, b \vee d) \in \theta(F)$. Similarly, using (GI1), one proves that $\theta(F)$ is compatible with $\wedge$.

Now, we prove that $\theta(F)$ is a gi-congruence. Let $(a,b), (c,d) \in \theta(F)$. Then there exist $f_{1}, f_{2} \in F$ such that $f_{1} \in (a \Rightarrow b) \cap (b \Rightarrow a)$ and $f_{2} \in (c \Rightarrow d) \cap (d \Rightarrow c)$. Let $x \in (a \Rightarrow c)$. Take $y=f_{1} \wedge f_{2} \wedge x$. Then, by (GI3) of Definition \ref{def: generalized implication},
\[
y \in (b \Rightarrow a) \cap (a \Rightarrow c) \subseteq b \Rightarrow c,
\]
i.e., $y \in (b \Rightarrow c)$. It follows that
\[
y \in (b \Rightarrow c) \cap (c \Rightarrow d) \subseteq b \Rightarrow d,
\]
i.e., $y \in b \Rightarrow d$. We take $f=f_{1} \wedge f_{2} \in F$. Since $y \leq x$, we have $y \Rightarrow x = A$, and $x \Rightarrow y = x \Rightarrow f$ by (GI1) and (GI4). Moreover, $1 \Rightarrow f \subseteq x \Rightarrow f$. Since $F$ is a gi-filter, $F \cap (1 \Rightarrow f) \neq \emptyset$, and so $F \cap (x \Rightarrow f) \neq \emptyset$. Thus, there exists $g \in F$ such that $g \in (x \Rightarrow y) \cap (y \Rightarrow x)$. Therefore, for each $x \in (a \Rightarrow c)$ there exist $y \in (b \Rightarrow d)$ and $g \in F$ such that $g \in (x \Rightarrow y) \cap (y \Rightarrow x)$, i.e., $(x,y) \in \theta(F)$, and $\theta(F)$ is a gi-congruence.

Let $a \in F$. Since $F$ is a gi-filter, $F \cap (1 \Rightarrow a) \neq \emptyset$, i.e., there is $f \in F$ such that $f \in (1 \Rightarrow a) = A \cap (1 \Rightarrow a) = (a \Rightarrow 1) \cap (1 \Rightarrow a)$. Then $(a,1) \in \theta(F)$ and $a \in [1]_{\theta(F)}$. So, $F \subseteq [1]_{\theta(F)}$.

Finally, let $\psi$ be a gi-congruence of $A$ such that $[1]_{\psi} \subseteq F$. Let $(a,b) \in \psi$. Since $(a,a), (b,b) \in \psi$ and $\psi$ is a gi-congruence, applying the definition to the pairs $(a,b), (a,a)$ we obtain that for $1 \in (a \Rightarrow a) = A$ there exists $x \in (b \Rightarrow a)$ such that $(1,x) \in \psi$, and applying it to the pairs $(b,a), (b,b)$ we obtain that for $1 \in (b \Rightarrow b) = A$ there exists $y \in (a \Rightarrow b)$ such that $(1,y) \in \psi$. Then $x,y \in [1]_{\psi} \subseteq F$ and $f=x \wedge y \in F$ satisfies $f \in (a \Rightarrow b) \cap (b \Rightarrow a)$, i.e., $(a,b) \in \theta(F)$ and $\psi \subseteq \theta(F)$.
\end{proof}

\begin{theorem}
Let $\langle A, \Rightarrow \rangle \in {\bf{GI}}$. If $\theta$ is a gi-congruence of $A$, then $\theta \subseteq \theta( [1]_{\theta} )$.
\end{theorem}
\begin{proof}
Since $[1]_{\theta}$ is a gi-filter, this follows from the last part of Theorem \ref{theo_congruence} with $\psi = \theta$ and $F = [1]_{\theta}$. For the sake of completeness, we give a direct proof. Let $(a,b) \in \theta$. Since $(a,a) \in \theta$ and $\theta$ is a gi-congruence, for $1 \in A = a \Rightarrow a$ there exists $x \in (b \Rightarrow a)$ such that $(1,x) \in \theta$, i.e., $x \in [1]_{\theta}$. Analogously, since $(b,b) \in \theta$, for $1 \in A = b \Rightarrow b$ there is $y \in (a \Rightarrow b)$ such that $(1,y) \in \theta$. So, $y \in [1]_\theta$. We take $z = x \wedge y \in [1]_{\theta}$. Then $z \in (a \Rightarrow b) \cap (b \Rightarrow a)$ and $(a,b) \in \theta([1]_{\theta})$. Therefore, $\theta \subseteq \theta([1]_{\theta})$.
\end{proof}

\begin{remark}
The inclusion $F \subseteq [1]_{\theta(F)}$ in Theorem \ref{theo_congruence} may be strict, as already happens in \cite[Example 4.5]{Ghorbani1}. Indeed, consider the commutative gi-lattice of Example \ref{gi_sd} and the gi-filter $F=[b)$. Since $b \in F \cap (0 \Rightarrow 1) \cap (1 \Rightarrow 0)$, we have $(0,1) \in \theta(F)$, and hence $[1]_{\theta(F)} = A \neq F$. Note, however, that $F$ is the kernel of a gi-congruence: the lattice congruence $\psi$ whose classes are $\{0,a\}$ and $\{b,1\}$ is a gi-congruence, because $x \Rightarrow y$ only depends on the $\psi$-classes of $x$ and $y$, and $[1]_{\psi} = F$.
\end{remark}

\begin{theorem}
Let $\langle A, \Rightarrow \rangle \in {\bf{CGI}}$. If $D$ is a gi-deductive system of $A$, then $D = [1]_{\theta(D)}$.
\end{theorem}
\begin{proof}
By Lemma \ref{gids implica gif}, $D$ is a gi-filter and by Theorem \ref{theo_congruence}, we know that $D \subseteq [1]_{\theta(D)}$. Let $a \in [1]_{\theta(D)}$. Then $(a,1) \in \theta(D)$ and there is $f \in D$ such that $f \in (a \Rightarrow 1) \cap (1 \Rightarrow a) = 1 \Rightarrow a$. So, $D \cap (1 \Rightarrow a) \neq \emptyset$ and, since $D$ is a gi-deductive system and $1 \in D$, $a \in D$. Therefore $D = [1]_{\theta(D)}$.
\end{proof}

If $\langle A, \Rightarrow \rangle \in {\bf{CGI}}$, then the map $D \mapsto \theta(D)$ from the set of gi-deductive systems of $A$ to the set of gi-congruences of $A$ is well-defined and injective, but in general it is not surjective.

\begin{example} \label{example_theta}
Let $A$ be the four-element chain $0<a<b<1$. The following table is the principal-ideal version of \cite[Example 4.10]{Ghorbani1}:
\vspace{0.4cm}
\begin{center}
\begin{tabular}{c|cccc}
$\Rightarrow$ & $0$ & $a$ & $b$ & $1$ \tabularnewline
\hline
$0$ & $A$      & $A$          & $A$   & $A$ \tabularnewline
$a$ & $\{0\}$  & $A$          & $A$   & $A$\tabularnewline
$b$ & $\{0\}$  & $A$          & $A$  & $A$\tabularnewline
$1$ & $\{0\}$  & $(a]$        & $(a]$  & $A$\tabularnewline
\end{tabular}
\par
\end{center}
\vspace{0.4cm}
Then $\langle A, \Rightarrow \rangle \in {\bf{CGI}}$ and $1 \Rightarrow c \subseteq (c]$, for all $c \in A$. By Lemma \ref{lemma_gif_gids}, the gi-filters and the gi-deductive systems coincide. Then the gi-deductive systems of $A$ are $\{1\}$, $[a)$ and $A$. On the other hand, the identity congruence $\Delta_{A} = \{(x,x) \colon x \in A\}$ is a gi-congruence of $A$. However, since $1 \in D \cap (a \Rightarrow b) \cap (b \Rightarrow a)$, we have $(a,b) \in \theta(D)$ for every gi-deductive system $D$, and so there is no gi-deductive system $D$ such that $\theta(D) = \Delta_{A}$.
\end{example}

\section{Priestley-style duality} \label{sec5}

We restrict the Priestley-style duality of \cite{Castro2} to commutative gi-lattices. The relational condition is the same as in the commutative WH-setting of \cite{Celani-Hernan}, but implication values are now ideals and the associated implication sets need only be open upsets. We begin by recalling gi-spaces.

\begin{definition}
A gi-space is a structure $\langle X, \leq, \tau, R \rangle$, where $\langle X, \leq, \tau \rangle$ is a Priestley space and $R$ is a binary relation on $X$ such that:
\begin{enumerate}
\item $R(x)$ is a closed set, for every $x \in X$.
\item $U \to_{R} V$ given by (\ref{impli conj}) is an open upset, for all $U,V \in CU(X)$.
\end{enumerate}
\end{definition}

The gi-lattice associated with a gi-space $\langle X, \leq, \tau, R \rangle$ is the pair $\langle CU(X), \Rightarrow_{R} \rangle$, where $\Rightarrow_{R}$ is given by $U \Rightarrow_{R} V = \{ W \in CU(X) \colon W \subseteq U \to_{R} V \}$, and $\to_{R}$ by (\ref{impli conj}).

\begin{lemma} \cite{Castro2} \label{aux_duality}
Let $\langle X, \leq, \tau, R \rangle$ be a gi-space. Then
\[
(x,y) \in R \Longleftrightarrow (\epsilon(x), \epsilon(y)) \in R_{CU(X)}.
\]
\end{lemma}

\begin{definition}
Extending the commutative WH-spaces of \cite{Celani-Hernan}, a commutative gi-space is a gi-space $\langle X, \leq, \tau, R \rangle$ that satisfies the following condition: for all $x,y,z \in X$, if $(x,y) \in R$ and $(y,z) \in R$, then there exist $u,v,w \in X$ such that $(x,u) \in R$, $(v,w) \in R$ and $y \leq w \leq z \leq u \leq v$.
\end{definition}

\begin{proposition} \label{prop_comm_space}
Let $\langle X, \leq, \tau, R \rangle$ be a gi-space. Then $\langle X, \leq, \tau, R \rangle$ is a commutative gi-space if and only if $\langle CU(X), \Rightarrow_{R} \rangle$ is a commutative gi-lattice.
\end{proposition}
\begin{proof}
Let $\langle X, \leq, \tau, R \rangle$ be a commutative gi-space. Suppose there exist $U,V,W \in CU(X)$ such that
\[
U \Rightarrow_{R} (V \Rightarrow_{R} W) \nsubseteq V \Rightarrow_{R} (U \Rightarrow_{R} W).
\]
Then there exists $S \in CU(X)$ such that $S \in U \Rightarrow_{R} (V \Rightarrow_{R} W)$ and $S \notin V \Rightarrow_{R} (U \Rightarrow_{R} W)$. Since $\langle CU(X), \Rightarrow_{R} \rangle$ is a gi-lattice, there is $P \in X(CU(X))$ such that $S \in P$ and $(V \Rightarrow_{R} (U \Rightarrow_{R} W)) \cap P = \emptyset$. By Lemma \ref{lema_F to I}, there is $Q \in X(CU(X))$ such that $V \in Q$, $(P,Q) \in R_{CU(X)}$ and $(U \Rightarrow_{R} W) \cap Q = \emptyset$. Again, by Lemma \ref{lema_F to I}, there is $T \in X(CU(X))$ such that $U \in T$, $(Q,T) \in R_{CU(X)}$ and $W \notin T$.
Then there exist $x,y,z \in X$ such that $\epsilon(x)=P$, $\epsilon(y)=Q$ and $\epsilon(z)=T$. So, $x \in S$, $y \in V$, $z \in U$ and $z \notin W$. It follows, by Lemma \ref{aux_duality}, that $(x,y) \in R$ and $(y,z) \in R$, and by hypothesis, there exist $u,v,w \in X$ such that $(x,u) \in R$, $(v,w) \in R$ and $y \leq w \leq z \leq u \leq v$.

On the other hand, by Lemma \ref{lema_F to I}, we have that $S \in U \Rightarrow_{R} (V \Rightarrow_{R} W) = [U) \Rightarrow_{R} (V \Rightarrow_{R} W)$. Then there are $F \in [U)$ and $I \in V \Rightarrow_{R} W$ such that $S \in F \Rightarrow_{R} I$, i.e., $U \subseteq F$, $I \subseteq V \to_{R} W$ and $S \subseteq F \to_{R} I$. Since $z \in U$, $z \leq u$ and $U$ is an upset, we have $u \in U \subseteq F$. Since $x \in S$, we have $x \in F \to_{R} I$ and $R(x) \cap F \subseteq I \subseteq V \to_{R} W$. Then $u \in R(x) \cap F$ and $u \in V \to_{R} W$, and since $V \to_{R} W$ is an upset, $v \in V \to_{R} W$, i.e., $R(v) \cap V \subseteq W$. As $y \in V$ and $y \leq w$, we have $w \in V$ and $w \in R(v) \cap V \subseteq W$. Hence, $w \in W$ and $w \leq z$ imply $z \in W$, which is a contradiction. Therefore, $U \Rightarrow_{R} (V \Rightarrow_{R} W) \subseteq V \Rightarrow_{R} (U \Rightarrow_{R} W)$. The other inclusion is similar and $\langle CU(X), \Rightarrow_{R} \rangle$ is a commutative gi-lattice.

Conversely, suppose that $\langle CU(X), \Rightarrow_{R} \rangle$ is a commutative gi-lattice. Let $x,y,z \in X$ be such that $(x,y) \in R$ and $(y,z) \in R$. By Lemma \ref{aux_duality} we have $(\epsilon(x), \epsilon(y)) \in R_{CU(X)}$ and $(\epsilon(y), \epsilon(z)) \in R_{CU(X)}$. Since $\langle CU(X), \Rightarrow_{R} \rangle$ is a commutative gi-lattice, by Theorem \ref{charact_2}, there exist $F,K,H \in X(CU(X))$ such that $(\epsilon(x),F) \in R_{CU(X)}$, $(H,K) \in R_{CU(X)}$ and $\epsilon(y) \subseteq K \subseteq \epsilon(z) \subseteq F \subseteq H$. It follows that there exist $u,v,w \in X$ such that $F=\epsilon(u)$, $K=\epsilon(w)$ and $H=\epsilon(v)$. Hence, $(\epsilon(x), \epsilon(u)) \in R_{CU(X)}$ and $(\epsilon(v), \epsilon(w)) \in R_{CU(X)}$ imply $(x,u) \in R$ and $(v,w) \in R$ by Lemma \ref{aux_duality}. Moreover, $\epsilon(y) \subseteq \epsilon(w) \subseteq \epsilon(z) \subseteq \epsilon(u) \subseteq \epsilon(v)$ and\bl{,} since $\epsilon$ is an order-isomorphism, we have $y \leq w \leq z \leq u \leq v$. Therefore, $\langle X, \leq, \tau, R \rangle$ is a commutative gi-space.
\end{proof}

Recall from \cite{Castro2} that if $\langle A, \Rightarrow_{A} \rangle, \langle B, \Rightarrow_{B} \rangle \in {\bf{GI}}$, then a homomorphism of bounded lattices $h \colon A \to B$ is a gi-homomorphism if
\[
{\rm{Idg}}(h[a \Rightarrow_{A} b]) = h(a) \Rightarrow_{B} h(b),
\]
for all $a,b \in A$. Composites of gi-homomorphisms are gi-homomorphisms and the class of gi-lattices with gi-homomorphisms forms a category. On the other hand, if $\langle X, \leq_{X}, \tau_{X}, R_{X} \rangle$ and $\langle Y, \leq_{Y}, \tau_{Y}, R_{Y} \rangle$ are two gi-spaces, then a function $f \colon X \to Y$ is a gi-morphism if it is a continuous and monotone function such that:
\begin{itemize}
\item If $(x,y) \in R_{X}$, then $(f(x),f(y)) \in R_{Y}$\bl{.}
\item If $(f(x),y) \in R_{Y}$, then there exists $z \in X$ such that $(x,z) \in R_{X}$ and $f(z)=y$.
\end{itemize}
Composites of gi-morphisms are gi-morphisms, and the category of gi-lattices with gi-homomorphisms is dually equivalent to the category of gi-spaces with gi-morphisms (see \cite{Castro2} for more details). Let ${\mathbb{CGI}}$ be the category whose objects are commutative gi-lattices and whose morphisms are gi-homomorphisms, and let ${\mathbb{CGIS}}$ be the category whose objects are commutative gi-spaces and whose morphisms are gi-morphisms. Then we have the following result.

\begin{theorem}
The categories ${\mathbb{CGI}}$ and ${\mathbb{CGIS}}$ are dually equivalent.
\end{theorem}
\begin{proof}
Since ${\mathbb{CGI}}$ and ${\mathbb{CGIS}}$ are full subcategories of the categories of gi-lattices and of gi-spaces, respectively, it suffices to check that the dual equivalence of \cite{Castro2} restricts to them. If $\langle A, \Rightarrow \rangle$ is a commutative gi-lattice, then its dual gi-space $\langle X(A), \subseteq, \tau_{A}, R_{A} \rangle$ is commutative by Theorem \ref{charact_2}. Conversely, if $\langle X, \leq, \tau, R \rangle$ is a commutative gi-space, then $\langle CU(X), \Rightarrow_{R} \rangle$ is a commutative gi-lattice by Proposition \ref{prop_comm_space}.
\end{proof}

{\bf{Declaration of Generative AI and AI-assisted technologies in the writing process:}} 
During the preparation of this manuscript, the authors used ChatGPT (OpenAI) and Claude (Anthropic) to assist with the organization of the manuscript, language editing and the preparation of portions of the LATEX source. All AI-assisted material was subsequently reviewed, verified and revised by the authors, who take full responsibility for the content of the paper.


\end{document}